\documentclass[leqno,12pt]{article}

\usepackage{geometry}
\usepackage{amssymb}
\usepackage{amsmath}

\newtheorem{Theorem}{Theorem}[section]
\newtheorem{Corollary}[Theorem]{Corollary}

\newtheorem{Proposition}[Theorem]{Proposition}
\newtheorem{Definition}[Theorem]{Definition}
\begin{document}
\title{Strong and safe Nash equilibrium in some repeated 3-player games }
\author{Tadeusz Kufel\\
\small Faculty of Economic Sciences and Management\\
\small N. Copernicus University in Toru\'n\\
\small e-mail: tadeusz.kufel@umk.pl
\and S\l awomir Plaskacz\\
\small Faculty of Mathematics and Computer Science\\
\small N. Copernicus University in Toru\'n\\
\small e-mail: plaskacz@mat.umk.pl
\and Joanna Zwierzchowska\\
\small Faculty of Mathematics and Computer Science\\
\small N. Copernicus University in Toru\'n\\
\small e-mail: joanna.zwierzchowska@mat.umk.pl}
\date{}
 \maketitle
{\bf JEL classification}\\ C72, C73, D91

 \vspace{3mm}

{\bf Key-words : }repeated game, strong Nash equilibrium, Blackwell's approachability, Lapunov function method.

 \vspace{3mm}

{\bf Abstract : }
The paper examines an infinitely repeated 3-player extension of the Prisoner's Dilemma game.
We consider a 3-player game in the normal form with incomplete information, in which each player has two actions. We assume that the game is symmetric and  repeated  infinitely many times. At each stage, players make their choices knowing only the  average payoffs from previous stages of all the players. A strategy of a player in the repeated game is a function defined on the convex hull of the set of payoffs. Our aim is to construct a strong Nash equilibrium in the repeated game, i.e.  a strategy profile being  resistant to deviations by coalitions. Constructed equilibrium strategies are safe, i.e. the non-deviating player payoff  is not smaller than the equilibrium payoff  in the stage game, and deviating players' payoffs do not exceed the non-deviating player payoff more than by a positive constant which can be  arbitrary small and chosen by the non-deviating player.  Our construction is inspired by Smale's good strategies described in \cite{smale}, where the repeated Prisoner's Dilemma was considered. In proofs we use arguments based on  approachability and strong approachability type results.

 \vspace{3mm}

\section{Introduction}\label{intro}

The Prisoner's Dilemma with its generalizations are very important as an example of conflicts and social dilemmas. As we can find in \cite{dawes}, social dilemmas are real life problems which have two properties: ''1. each individual receives a higher payoff for a socially defection choice than for a socially cooperative choice, no matter what the other individuals in society do; 2. all individuals are better off if all cooperate than if all defect.'' An example of such situation in real life is a problem of soldiers who fight in a battle. They are personally better off taking no chances, yet if no one fight against the enemy, then the result will be worst for all of soldiers. Such dilemmas can be found among resource depletion, pollution and overpopulation.

Social dilemmas are games in which there is a conflict between individual rationality and optimality of the equilibrium payoff. Since it is observable that people cooperate with each other in the real situations, game theorists have faced the obstacle, how to construct simple tools to encourage players in such games to cooperate with each other. The model need to approximate the real situation and strategies should be likely to use.

A natural approach is to consider the infinitely repeated game. Usually, all players observe the whole history of action profiles used in previous stages of the repeated game. Such situation is called the game with complete information. The strategies are functions from the set of the histories into the set of actions. Payoffs in the repeated game are either the discounted sum of stage payoffs or the limit of average payoffs. The aim of this approach is to obtain the Nash equilibrium in the repeated game with the pair of payoffs which is  close to the cooperation payoffs in the stage game. Since the fifties of the last century there appeared various \textit{folk theorems} which was not explicitly published and, in many cases, the original author is unknown.

The classic Prisoner's Dilemma is a 2-player game, in which each player has two actions, usually denoted as $C$ (cooperation) and $D$ (defection). The game has a unique Nash equilibrium -- a pair of actions such that the action of each of the players optimize this player's payoff given the action of the opponent. The Nash equilibrium is the action profile $(D,D)$ which is the pair of strictly dominant actions i.e. playing $D$ is better than $C$ whatever the other player does. What is more, both players benefit changing $(D,D)$ into $(C,C)$. So, the mechanism of individual rationality fails in the Prisoner's Dilemma and it leads to a loss of both players. It means that the Nash equilibrium is not Pareto-optimal in this case.

One of solutions for lack of cooperation of the Nash equilibrium in the stage game is an idea of good strategies introduced by Smale in \cite {smale} for the repeated Prisoner's Dilemma. Every pair of good strategies is a Nash equilibrium in the repeated game with Pareto-optimal payoffs corresponding to the payoff of $(C,C)$ in the stage game. The second advantage of the good strategies equilibrium is the warranted minimal payoff for the non-deviating player. The minimal payoff is equal to the Nash payoff in the stage game. Good strategies have yet another advantage that has not been pointed in \cite{smale}. Choosing a good strategy appropriately, the player  controls  the second player's payoff. For every $\varepsilon>0$ there exists the $\varepsilon$-good strategy of the first player such that for an arbitrary second player's strategy, the first player's payoff will be at most $\varepsilon$ smaller than the second player's payoff. The Prisoner's Dilemma is symmetric, so the second player also can choose the $\varepsilon$-good strategy which provides him no worse payoffs than the first player's one minus  $\varepsilon$.

In fact, good strategies have properties postulated by R. Axelrode in \cite{ax}. In 80's Axelrode studied the evolution of cooperation. It refers to how cooperation can emerge and persist as elucidated by application of game theory. He organized a tournament in which game theory experts submitted their strategies and each strategy was paired with each other for 200 iterations of Prisoner's Dilemma. Accumulated payoffs through the tournament was treated as a score. The winner was the strategy submitted by Arnold Rappaport -- Tit for Tat. The additional advantage of this tournament was detecting what properties strategies should satisfy to encourage players to cooperate. They should be: nice, forgiving, retaliatory and are founded on simple rules. Good strategies have these properties and, what is more, player cooperates until the other player's average payoff is greater than his average payoff plus $\varepsilon$. By choosing  $\varepsilon$, the player determines  the level of his tolerance for the defection.

In this paper we shall consider the generalization of the idea of good strategies onto the Prisoner's Dilemma type repeated game for three players.
We consider the repeated game with a partial monitoring. We assume that, after each stage, all players can only observe  an aggregated history -- the arithmetic mean of the payoffs from previous stages. The stage game is a symmetric 3-player game where each player has an action set  consisting of two actions : $I$ and $NI$\footnote{From now on, we choose to name strategies with $I$ and $NI$, where $I$ means invest and it corresponds to strategy $C$ and $NI$ corresponds to strategy $D$. }.
We assume that the action  profile $(NI,\,NI,\,NI)$ is the only Nash equilibrium, and the sum of the players payoffs is minimal for this profile. The sum of players payoffs is maximal for the profile $(I,\,I,\,I)$. The strategy profile in the repeated game is a function $s:S\to\{I,\,NI\}^3$, where $S\subset R^3$ is a convex hull of the set of the payoffs in the stage game.

Our aim is to construct a strategy profile $s^*$ which is an approximated  strong Nash equilibrium in the repeated game under consideration. The constructed equilibrium is safe in the meaning that the payoff of a player choosing strategy $s^*_i$ is not less then the equilibrium payoff in the stage game. This payoff is assured even if the other two players choose an arbitrary strategy. Furthermore, the $\varepsilon$-good strategy guarantees that, in long time horizon, other player's average payoff will not exceed the good strategy player's average payoff by more than $\varepsilon$.

 The notion of the strong equilibrium in the framework of repeated games was introduced by Aumann (see \cite{Aum}, \cite{Aum2}, \cite{Aum3}), who showed that every payoff that belongs to the $\beta$-core of the stage game is a strong equilibrium payoff in the corresponding repeated game (see \cite[Thm. 6.2.2]{Sor1} ). Despite the fact that the payoff corresponding to the profile $(I,\,I,\,I)$ belongs to the $\beta$-core, our result is not exactly a case of the Aumann results. We have dropped the assumption of the full monitoring. Players do not observe the full history, i.e. the sequence of actions selected by all players in the previous periods. Instead, we assume that  they observe the aggregate history, i.e. the arithmetic mean of the previous payoffs of all the players.  It is worth noting that the results on  the existence of strong equilibria (see \cite{kbw}, \cite{nt}) do not apply to the repeated game considered in the present paper.

 The repeated Prisoner's Dilemma for more than two players has been considered in \cite{bbh}. The $\varepsilon$-good strategies constructed in the paper have some additional properties to the strategies in \cite{bbh}.  In \cite{bbh} the authors base on similar approchability results as we do in this paper. The difference is that authors consider $N$-players Prisoner's Dilemma Game in which strategies are stochastic processes. In our approach all strategies are deterministic.

The paper is organized as follows.

In Section \ref{sec:1} we present the basic information about sequences  related to a map of a convex set. We adopt Blackwell's approachability method (see \cite{black}) which was originally  used in the framework of 2-player repeated games with vector payoffs. We show that the Blackwell condition is sufficient to obtain the convergence of the sequence of arithmetic means to a set  called a weak attractor. The weak attractors introduced in Subsection \ref{sec:1.1} have different properties in comparison with approachable sets in the sense of Blackwell. We provide an example of a singleton being a weak attractor that does not satisfy the Blackwell condition. Such a situation is not possible for approachable sets (comp. \cite[Thm. 8]{Shani}).
In repeated games, there is considered a sequence of vector payoffs. Each payoff corresponds to one repetition of the state game.  Subsection \ref{sec:1.1} provides us necessary results to analyze the directions in which the trajectory shifts and to examine the convergence of such sequence. This is crucial for defining the payoff in the repeated game.
Subsection \ref{sec:1.2} provides basic properties of the Banach limit which shall be used to prove that $\epsilon$-good strategies are $\epsilon$ Nash equilibria.
 In some  of our arguments we not only require   that the sequence of mean payoffs converges  to a set, but that almost all its entries belong to the set.  A similar problem named strong approachability was considered in \cite{Shani}. In Section \ref{sec:2} we adopt a  Lyapunov  function method for discrete and discontinuous dynamical systems to obtain a deterministic strong approachability result.

In Section \ref{sec:3} we consider a repeated 3-player symmetric game. Every player has two actions: invest ($I$) or not invest ($NI$). The vector payoff $B=(p_3,\,p_3,\,p_3)$ corresponding to the strategy profile $(I,I,I)$ is Pareto optimal and the strategy profile $(NI,NI,NI)$ is a Nash equilibrium in the stage game with the payoff vector $(r_0,\,r_0,\,r_0))$. We assume that,  in the repeated game, every player knows the average vector payoff from the previous stages of the game. The strategy $s_i:S\to\{I,NI\}$, $i\in\{1,2,3\}$, is a function from the convex hull $S$ of vector payoffs in the stage game to the set of his actions $\{I,NI\}$.
The strategy profile $s=(s_1,s_2,s_3)$ and the vector payoff function $G:\{I,NI\}^3\to R^3$ determine the function $\varphi=G\circ s:S\to S$. The strategy profile  $s$ and the initial point $x_1\in S$ determine the trajectory $\bar{x}_n(s,x_1)$  of a dynamic system given by
\[
\bar{x}_{n+1}=\frac{n\bar{x}_n+\varphi(\bar{x}_n)}{n+1}
\]

Our aim is to construct a strategy profile $s_\varepsilon^*=(s_1^*,\,s_2^*,s_3^*)$ such that for every $x_1\in S$
\begin{equation}\label{war1}
\lim_{n\to\infty}\bar{x}_n=B,
\end{equation}
where $\bar{x}_n=\bar{x}_n(s_\varepsilon^*,x_1)$. If one player (for example player 3) deviates then
\begin{equation}\label{war2}
 \limsup_{n\to\infty}\bar{x}_n^3\leq p_3+\varepsilon,
\end{equation}
where $\bar{x}_n=\bar{x}_n((s_1^*,s_2^*,s_3),x_1)$ and $s_3:S\to\{I,NI\}$ is an arbitrary strategy of player 3. If two players deviate (for example players 2 and 3) then
\begin{equation}\label{war3}
\limsup_{n\to\infty}(\bar{x}_n^2+\bar{x}_n^3)\leq 2 p_3,
 \end{equation}
 \begin{equation}\label{war4}
 \liminf_{n\to\infty}\bar{x}_n^1\geq r_0
 \end{equation}
 \begin{equation}\label{war5}
 \lim_{n\to\infty} \mbox{dist}(\bar{x}_n,\{x\in S;\,x_2\leq x_1+\varepsilon,\;x_3\leq x_1+\varepsilon\})=0
\end{equation}
where $\bar{x}_n=\bar{x}_n((s_1^*,s_2,s_3),x_1)$ and $s_2,\,s_3:S\to\{I,NI\}$ are the arbitrary strategies of players 2 and 3 respectively.\\
If the payoff is a Banach limit (comp.\cite{conway}) of the sequence of average payoffs then the strategy profile $s_\varepsilon^*$ is a strong $\varepsilon$-Nash equilibrium in the repeated game as a consequence of (\ref{war1}-\ref{war3}). Property (\ref{war4}) implies that the non-deviating player's payoff is no smaller than the payoff corresponding to the Nash equilibrium in the stage game. Property (\ref{war5}) guarantees that the deviating player's payoff will not exceed the good strategy player's payoff by more than $\varepsilon$.  The results presented in Theorems \ref{T3}, \ref{T4}, \ref{T2} give a partial answer to the question asked by Smale in the last Remark in section 1 of \cite[p. 1623]{smale}.

Section \ref{sec:4} contains concluding remarks. In Appendix \ref{sec:A} we present the proofs of theorems from Subsection \ref{sec:2}.

\section{Preliminaries}\label{sec:1}

\subsection{Approachability  results}\label{sec:1.1}

 Let $H$ be a finite
dimensional vector space and $\langle\cdot,\cdot\rangle$,
$|\cdot|$ denote an inner product and a norm in $H$, respectively. We assume
that  $S$ is a nonempty convex closed subset of $H$.
 By $N^\varepsilon(B)$  we denote an $\varepsilon$-neighbourhood of the
set $B$ in $S$, i.e. $N^\varepsilon(B)=\{x\in
S:\,dist(x,B)<\varepsilon\}$.  The closure (the convex hull) of the
set $A$ we denote by $cl(A)$ ($co(A)$) .

We study limit properties of sequences
$(\bar{x}_n)_{n=1}^\infty$ defined by a map $\varphi:S\to S$ and
an initial point $x_1\in S$ by
\begin{equation}\label{traj}
\bar{x}_{n+1}=\frac{n\bar{x}_n+\varphi(\bar{x}_n)}{n+1},\;\;\bar{x}_1=x_1,
\end{equation}
  The sequence $(\bar{x}_n)$ can be interpreted as a sequence of arithmetic means
 $\bar{x}_n=\frac{1}{n}(x_1+\ldots+x_n)$, where $x_{k+1}=\varphi(\bar{x}_k)$. The map $\varphi$ defines a dynamical system $\beta_n:S\to S$ by
 \[
 \beta_n(x)=\frac{n x+\varphi(x)}{n+1},\; n=1,\,2,\,\ldots
 \]
 We denote by $\bar{x}_n(\varphi,\,x_1)$ a trajectory determined by (\ref{traj}).

We say that a closed set $A\subset S$ is a \textit{weak attractor} for a dynamic system determined by the map $\varphi$ if for every $x_1\in S$ we have
\[
\lim_{n\to\infty} dist(\bar{x}_n(\varphi,\,x_1),\,A)=0
\]
where $dist(\cdot,A)$ denotes the distance to the set $A$. We
provide some sufficient conditions  for being a weak attractor.

First we
formulate   Blackwell approachability type theorem that originally
was presented in  \cite{black}  in the framework of  repeated
games with vector payoffs. We say that a map $\varphi:S\to S$
\textit{satisfies the Blackwell condition for a set $A\subset S$ in the domain
$D\subset S$} if
\begin{equation}\label{condition B}
\forall x\in D,\,\exists y\in \Pi_A(x),\;\;\langle x-y,\varphi(x)-y\rangle\leq 0,
\end{equation}
where $\Pi_A(x)$ denote the set of points in $A$ that are proximal to $x$, i.e. $\Pi_A(x)=\{a\in A:\,|a-x|=dist(x,A)\}$.

  The  deterministic version of  the Blackwell
approachability result can be formulated in the following way.

\begin{Proposition}\label{bt}
Suppose that the  map $\varphi:S\to S$ satisfies the Blackwell condition for a  closed set $A\subset S$ in the domain $D\subset S$. If almost all elements of the  bounded sequence $\bar{x}_n(\varphi,x_1)$ belong do the set D then
\[
\lim_{n\to\infty} dist (\bar{x}_n,A)=0
\]

\end{Proposition}
We provide the proof of Proposition \ref{bt} for  the reader convenience.

{\bf Proof : }
For a sufficiently large $n$  we choose $y\in \Pi_K(\bar{x}_n)$ and then
\[
dist(\bar{x}_{n+1},A)^2\leq
|\bar{x}_{n+1}-y|^2=\left|\frac{n}{n+1}(\bar{x}_n-y)+\frac{1}{n+1}(\varphi(x_{n})-y)\right|^2=
\]
\[
=\left(\frac{n}{n+1}\right)^2|\bar{x}_n-y|^2+\left(\frac{1}{n+1}\right)^2|\varphi(x_{n})-y|^2+2\frac{n}{(n+1)^2}
\langle\bar{x}_n-y,\varphi(x_{n})-y\rangle
\leq
\]
\[
\leq\left(\frac{n}{n+1}\right)^2 dist(\bar{x}_n,A)^2+\left(\frac{1}{n+1}\right)^2 C
\]
where $C$ is an upper bound of $dist(x_n,A)$. Setting $d_n=n^2
dist(\bar{x}_n,A)^2$ we have $d_{n+1}\leq d_n+C$ for $n\geq n_0$. Thus $d_n\leq
d_{n_0}+(n-n_0)C$. So
\[
dist(\bar{x}_n,A)^2\leq\frac{1}{n}\left(d_{n_0}+\frac{n-n_0}{n} C\right).
\]
\begin{flushright} QED \end{flushright}

\begin{Corollary}\label{2bc}
If the map $\varphi:S\to S$ satisfies the Blackwell condition for a closed set $A\subset S$ in the domain $S$, then the set $A$ is a weak attractor for $\varphi$. If the set $A\subset S$ is convex and the map $\varphi:S\to A$ maps into the set $A$ then the set $A$ is a weak attractor for $\varphi$.
\end{Corollary}

Taking $A=(-\infty,c]$ in Proposition \ref{bt} we obtain the
following property of real sequences.
\begin{Corollary}\label{lBC}
Suppose that   $(a_n)_{n=1}^\infty$ is a bounded sequence in
 $ \mathbb{R}$ and $(\bar{a}_n)_{n=1}^\infty$ is the sequence of arithmetic means, i.e. $\bar{a}_n=\frac{1}{n}\sum_{k=1}^n a_k$.If we have

 \[
(\bar{a}_n>c \quad \Rightarrow \quad a_{n+1}\leq c)
 \]
for almost all $n$   and a fixed constant $c\in \mathbb{R}$, then
 \[
 \limsup_{n\rightarrow \infty }\bar{a}_n \leq c.
 \]
 \end{Corollary}

In many cases the set $A$ is a weak attractor despite that the
Blackwell condition is not satisfied. Such a situation occurs in repeated games that we study in Section \ref{sec:3}. Below we present two properties of weak attractors which are necessary for our reasoning.

\begin{Proposition}\label{atractor1}
Suppose that the sets $A,\,B\subset \mathbb{R}^d$ are nonempty  closed and $B$ is bounded. If a sequence $x_n$ satisfies
\[
\lim_{n\to\infty} dist (x_n,A)=\lim_{n\to\infty} dist (x_n,B)=0
\]
then
\[
\lim_{n\to\infty} dist (x_n,A\cap B)=0
\]
\end{Proposition}
{\bf Proof : }
We choose $a_n\in A$, $b_n\in B$ such that
\[
|x_n-a_n|=dist(x_n,A),\;\;|x_n-b_n|=dist(x_n,B)
\]
Since the set $B$ is compact, we obtain that the sequences $(a_n)$, $(b_n)$, $(x_n)$ are bounded and they have the same nonempty set
$C$ of accumulating points. Thus $\lim_{n\to\infty} dist (x_n,C)=0$ and $C\subset A\cap B$.
\begin{flushright} QED \end{flushright}

\begin{Proposition}\label{atractor2}
We suppose that a closed set $A\subset S$ is a weak attractor for
the map $\varphi:S\to S$ and a closed subset $B\subset A$ satisfies
\begin{equation}\label{ed}\begin{array}{ll}
\forall \varepsilon>0,\,\exists \delta>0,\;&\mbox{ $\varphi$
satisfies the Blackwell condition}\\
&\mbox{ for the set $cl(N^\varepsilon(B))\cap A$ in the domain
$N^\delta(A)$ .}
\end{array}
\end{equation}
Then the set $B$ is a weak attractor for $\varphi$.
\end{Proposition}
{\bf Proof : }
Fix $x_1\in S$ and $\varepsilon>0$. By (\ref{ed}), we choose
$\delta>0$ such that almost all elements of the trajectory
$\bar{x}_n(\varphi,\,x_1)$ belongs to $N^\delta(A)$. By
Proposition \ref{bt}, we obtain
\[
\lim_{n\to\infty} dist(\bar{x}_n,cl(N^\varepsilon(B))\cap A)=0
\]
Thus
\[
\limsup_{n\to\infty} dist(\bar{x}_n,B)\leq\varepsilon
\]

\begin{flushright} QED \end{flushright}

The method illustrated in Proposition \ref{atractor1} and Proposition \ref{atractor2} bases on the scheme that we explain  in the following
example.\\

{\bf Example 1}
Let $S=\mathbb{R}^2$, $a,b\in \mathbb{R}^2$, $a_2<0$, $b_2>0$, $a_1\neq b_1$ and
\[
\varphi(x,y)=
\left\{ \begin{array}{cc}
a&\mbox{ if } y>0,\\
b&\mbox{ if } y\leq 0.
\end{array}\right.
\]
We show that $\lim_{n\to\infty}\bar{x}_n=d$ for every $\bar{x}_1\in \mathbb{R}^2$, where the limit $d$ is the point of intersection of the interval $\overline{ab}$ with the line $p=\{(x,y):\,y=0\}$. The set $D=\{d\}$ does not satisfy condition (\ref{condition B}). Indeed, if $a_1<b_1$ and $x>d_1$ then $\varphi(x,0)=b$ and $\langle (x,0)-(d_1,d_2),\varphi(x,0)-(d_1,d_2)\rangle> 0$. To show that the set $D$ is a weak attractor we point out  weak attractors $A$, $B$  such that $D=A\cap B$. We set $A=p$ and $B=\overline{ab}$. The sets $A$, $B$ satisfy the Blackwell condition (\ref{condition B}). By Theorem  \ref{bt}, we have
\[
\lim_{n\to\infty} dist (\bar{x}_n,A)=\lim_{n\to\infty} dist (\bar{x}_n,B)=0
\]
Applying Proposition \ref{atractor1} we obtain that $\lim_{n\to\infty}\bar{x}_n=d$.\\

Finally, we shall formulate a property of the dynamical system.

\begin{Proposition}\label{l1}
If the set $S$ is bounded then for every $\xi >0$  there exists $N\in \mathbb{N}$ such that for
all $n>N$ and for all $x\in S$
\begin{displaymath}
|\beta_n(x)-x|<\xi,
\end{displaymath}
where the map $\varphi:S\to S$ determining $\beta_n$ is arbitrary.
\end{Proposition}

\subsection{Payoff in the repeated game}\label{sec:1.2}

Considering a sequence of payoffs in the repeated games we always receive a bounded sequence. As we presented in $(\ref{traj})$, the dynamic is the vector of the arithmetic mean of the payoffs received in the previous repetitions. To analyze such sequence, the following proposition shall be useful.

\begin{Proposition}\label{lA}
Suppose that $a_0,a_1,\dots, a_k \in \mathbb{R}^d$ and let $T\in
\mathbb{N}$. Then for all $\epsilon >0$  and for all $n_1,\dots,
n_k \geq 0 $ such that $n_1+\dots+n_k=n$, where $n$ is
sufficiently large,  we have
\begin{displaymath}
\frac{T}{T+n}a_0+\frac{n_1}{T+n}a_1+\dots + \frac{n_k}{T+n}a_k \in
N_{\epsilon}(co\{a_1,\dots, a_k\}).
\end{displaymath}
\end{Proposition}
Proposition \ref{lA} is a consequence of the fact that $\frac{T}{T+n}a_1+\frac{n_1}{T+n}a_1+\dots + \frac{n_k}{T+n}a_k \in
co\{a_1,\dots, a_k\}$ and $\frac{T}{T+n}|a_0-a_1|$ is small where $n$ is
sufficiently large.

To define the payoff in repeated games we shall use  the Banach limit (see \cite{conway}).  The Banach limit $L$ is a continuous linear functional definite on the space $l^\infty$ of bounded scalar sequences. If $(x_n)$ is a bounded sequence of points in $R^d$ then $\mbox{Lim}(x_n):=(\mbox{Lim}(x_{n1}),\mbox{Lim}(x_{n_2}),\ldots, \mbox{Lim}(x_{nd}))$, where $x_n=(x_{n1}, x_{n2},\ldots,x_{nd})$. So  Banach Limit can be extended onto the space of bounded sequences of points in $R^d$. If $\varphi:R^d\to R$ is a linear functional then $\varphi(\mbox{Lim}(x_n))=\mbox{Lim}(\varphi(x_n))$.

\begin{Proposition}\label{banachlimit}
If $A$ is a compact convex subset of $R^d$ and a sequence $(x_n)\subset R^d$ satisfies $\lim_{n\to\infty} dist(x_n,A)=0$,
then $\mbox{Lim}(x_n)\in A$.
\end{Proposition}

{\bf Proof.} Suppose to the contrary that $\mbox{Lim}(x_n)\notin A$. Then there exists a functional $\varphi:R^d\to R$ such that $\varphi(\mbox{Lim}(x_n))>\sup_{a\in A}\varphi(a)$. We have $\limsup_{n\to\infty} \varphi(x_n)\leq\sup_{a\in A}\varphi(a)$. Thus
\[
\varphi(\mbox{Lim}(x_n))=\mbox{Lim}(\varphi(x_n))\leq \limsup_{n\to\infty} \varphi(x_n)\leq \sup_{a\in A}\varphi(a)
\]
which gives the contradiction.

\begin{flushright} QED \end{flushright}

\subsection{A Lapunov type result}\label{sec:2}
The Lapunov function method is typically used to study stability
of equilibrium points  for dynamical systems. Using the Lapunov function method we obtain a strong approachability result for
a dynamical system determined by a multivalued map.

Let $H$ be a Hilbert space and $p_1,...,p_k\in H$ be unit vectors, i.e. $|p_i|=1$. We define
a function $V:H\to \mathbb{R}$ by
\begin{equation}\label{function V}
V(x)=\max_{i\in\{1,\,\ldots,\,k\}}V_i(x) \qquad \mbox{where}
\qquad V_i(x)=\langle p_i,x\rangle.
\end{equation}
The function $V$ is a support function of the set
$\{v_1,\dots,v_k\}$. So, the function $V$ is convex, positively homogeneous and
lipschitz continuous  with the
 constant $L=1$ (see \cite{rock}).

Set
\begin{displaymath}
\Delta_c=\bigcap_{i=1}^k \{x\in H:V_i(x)<c\}=\{x\in H:V(x)<c\}.
\end{displaymath}
 Let us denote by $\varphi:S\multimap S$ a
multivalued map of a subset $S\subset H$.

\begin{Definition}\label{DLapunov}
 We say that $V$ is the Lapunov type function for the multivalued map
$\varphi$ with the constant $c>0$ if
\begin{equation}\label{Lapunov}
\exists 0<\delta <c ,\, \forall x\in S\setminus \Delta_c ,\,
\forall i=1,...,k ,\,\forall \omega \in \varphi(x) ,\;
(V_i(x)\geq V(x) - \delta \; \Rightarrow  \; V_i(\omega) \leq 0)
\end{equation}
\end{Definition}

If the function $V$ satisfies
\begin{equation}\label{Lapunov1}
\forall x\in S\setminus\Delta_c,\,\forall
\omega\in\varphi(x),\,\forall i\in\{1,\ldots, k\}\;\;V_i(x)>0\; \Rightarrow\;V_i(\omega)\leq 0
\end{equation}
then $V$ is  the Lapunov type function for $\varphi$ with the
constant $c$.

If $V$ is the Lapunov type function for $\varphi$ with the
constant $c$ and $c_1>c$ then $V$ is the Lapunov type function
for $\varphi$ with the constant $c_1$.

To explain why we say that $V$ is the Lapunov type function observe that if $V_i(x)=V(x)$ then $p_i\in\partial V(x)$, where $\partial V(x)$
is the subdifferential of a convex function. The condition
(\ref{Lapunov}) implies the following inequality
\[
\langle p_i , \omega - x\rangle \leq \langle p_i , \omega
\rangle - \langle p_i , x\rangle \leq 0-V(x)+\delta <\delta-c
 < 0 \quad \mbox{for } \quad x\in S\setminus \Delta_c
\]
which  means that $V$ is the Lapunov function for the vector field
$f(x)=\omega - x$.

\begin{Proposition}\label{lfm}
Let $S$  be a nonempty bounded convex subset of $H$ and the
function $V:H\to\mathbb{R}$ given by (\ref{function V}) be the Lapunov
type function for the multivalued map $\varphi:S\multimap S$ with
the constant $c>0$. If  a sequence $(\bar{x}_n)_{n=1}^\infty$
satisfies
\begin{equation}\label{zm}
\bar{x}_1=x_1\in S,\;\;\bar{x}_{n+1}=\frac{n\bar{x}_n+
x_{n+1}}{n+1},\;\;x_{n+1}\in\varphi(\bar{x}_n)
\end{equation}
then
\[
\forall c_1>c,\,\exists N,\; \forall n\geq N,\;\;
\bar{x}_n\in\Delta_{c_1}
\]
\end{Proposition}

The proof of Proposition \ref{lfm} is technical and it is presented in Appendix \ref{sec:A}.

\section{The model and main results}\label{sec:3}
Let $G$ be a 3-player symmetric game and every player has two pure actions: "invest" ($I$) or "not invest" ($NI$). By $P_I$ ($P_{NI}$)
we denote the payoff for an investing (not investing)  player.  All payoffs depend
on the total number of investing players. If $n\in\{0,\,1,\,2,\,3\}$
is the total number of investing players, then
\begin{displaymath}
\begin{array}{c|c|c}
n & P_I(n) & P_{NI}(n) \\ \hline
0 & - & r_0 \\
1 & p_1 & r_1 \\
2 & p_2 & r_2 \\
3 & p_3 & -
\end{array},
\end{displaymath}
The game $G$ in the normal form is given by the matrix:
\begin{displaymath}
\begin{array}{c|cc}
I & (p_2,r_2,p_2) & (p_3,p_3,p_3) \\
NI & (r_1,r_1,p_1) & (r_2,p_2,p_2) \\
\hline
 & NI & I
\end{array}.
\end{displaymath}
when the third player invests, and by the matrix
\begin{displaymath}
\begin{array}{c|cc}
I & (p_1,r_1,r_1) & (p_2,p_2,r_2) \\
NI & (r_0,r_0,r_0) & (r_1,p_1,r_1) \\
\hline
 & NI & I
\end{array},
\end{displaymath}
when the third player does not invest.

We shall assume that  the functions $P_I(\cdot)$, $P_{NI}(\cdot)$ are increasing:
\begin{equation}
0<r_0<r_1<r_2 \qquad \mbox{ and } \qquad 0<p_1<p_2<p_3. \label{nierpodst}
\end{equation}
  We assume that
\begin{equation}
p_1<r_0. \label{nierNash}
\end{equation}
By (\ref{nierNash}),   the outcome $(NI,NI,NI)$ is a Nash equilibrium.
We assume that the more players invest, the greater the sum of all players payoffs is, i.e.
\begin{equation}
3r_0<p_1+2r_1<2p_2+r_2<3p_3.\label{nierinvest}
\end{equation}
By (\ref{nierinvest}), the vector payoff $(p_3,p_3,p_3)$ is  Pareto optimal. In fact, the condition $(\ref{nierinvest})$ means even more -- the vector payoff  $(p_3,p_3,p_3)$ maximize the sum of payoffs. To obtain a strong equilibrium in the repeated game we assume that:
\begin{equation}
p_1+r_1<2p_3. \label{nierstrong}
\end{equation}
We additionally assume  that:
\begin{equation}
p_2<r_2. \label{nierantypod}
\end{equation}
Observe that from the opposite inequality  $r_2\leq p_2$ implies that   $(p_3,p_3,p_3)$  is a Nash equilibrium payoff, what we wanted to avoid.

We introduce the following notations
\begin{displaymath}
\begin{array}{c}
A=(r_0,r_0,r_0),\\
B=(p_3,p_3,p_3), \\
C_1^1:=(p_1,r_1,r_1),\\
C_2^1:=(r_1,p_1,r_1),\\
C_3^1:=(r_1,r_1,p_1),\\
C_1^2:=(r_2,p_2,p_2),\\
C_2^2:=(p_2,r_2,p_2),\\
C_3^2:=(p_2,p_2,r_2).
\end{array}
\end{displaymath}
If $i$ players invest ($i\in\{1,\,2\}$) then  $C^i_j$ denotes the vector payoff in the game $G$. If
$i=1$ then $j$ shows which one invests, while if $i=2$ then $j$ tells which player does not invest.

The strategy profile in the iterated game is given by a map $s:S\to\{I,\,NI\}^3$, where $S$ is the convex hull of vector payoffs set, i.e.
\[
S=conv\{A,\,B,\,C^1_1,\,C^1_2,\,C^1_3,\,C^2_1,\,C^2_2,\,C^2_3\}
\]
 The strategy profile $s$ determines a dynamical process $\beta_n:S\rightarrow S$
\begin{equation}
\beta_n(x)=\frac{nx+\varphi(x)}{n+1}, \qquad \mbox{ for } \quad x\in S,\quad
T\in \mathbb{N}, \label{dynamika}
\end{equation}
where $\varphi:S\rightarrow S$ is given by the formula $\varphi = G\circ s$. Observe that a pair $(s,\,x_1)$, where $s$ is a strategy profile and $x_1\in S$, uniquely determines  a sequence $(\bar{x}_n)_{n=1}^\infty$ by:
\[
\bar{x}_1=x_1,\;\;\;\bar{x}_{n+1}=\beta_n(\bar{x}_n)
\]
We denote the obtained sequence by $\bar{x}_n(s,x_1)$. A similar construction of a sequence was considered in Section \ref{sec:1}. The strategy profile $s$ and the initial point $x_1\in S$ uniquely determine a play path. The action profile in the next stage $s(\bar{x}_n)$ depends on the average  vector payoff  $\bar{x}_n$. The element  $x_{n+1}$ is the vector payoff in $n+1$ stage. We do not assume that the players observe the full history of the game. Instead, they observe aggregated history -- the arithmetic mean of vector payoffs.

Motivated by the Smale construction in \cite{smale} we define an $\epsilon$-good
strategy for the $i$-th player $s_i^\varepsilon :S\to\{I,\,NI\}$  by
\begin{equation}\label{gs}
s_i^\varepsilon(x)=\left\{
\begin{array}{ll}
I & \mbox{ if } x\in V_i \\
NI & \mbox{ if } x\in S\setminus V_i
\end{array}
\right.
\end{equation}
where
\begin{displaymath}
\begin{array}{cl}
V_i=&\Omega_i^{\epsilon}\setminus W_i,\\
\Omega_i^{\epsilon}=&\{x\in S:\quad x_i>x_j-\epsilon\mbox{ and } x_i>x_k-\epsilon \},\\
W_i=&\{x\in S:\quad x_i<r_0 \mbox{ or } x_j+x_k>2p_3 \}
\end{array}
\end{displaymath}
where $i$, $j$, $k$ are pairwise different elements of the set of players $\{1,\,2,\,3\}$.
The player  invests if his average payoff is greater than the every other players' average payoff  minus $\varepsilon$. The player stops investing if his playing $I$ has been exploited by his opponents, that is either the average payoff of the player is lower than the payoff guaranteed by Nash equilibrium  ($x_i<r_0$) or the sum of the other players' average payoffs  is greater then the sum of their payoffs corresponding to the Pareto optimal profile $(I,I,I)$ ($x_j+x_k>2p_3)$.

First we consider the case when all  players choose good strategies.
Then the average payoffs vector tends to the point $B$ corresponding to the  Pareto optimal profile $(I,I,I)$.
\begin{Theorem}\label{T3}
Suppose that $s_i^{\varepsilon}:S\rightarrow \{I,NI\}$ are the $\epsilon$-good
strategies for $i=1,2,3$. Then
\[
\lim_{T\to\infty}\bar{x}_T=B,
\]
where $\bar{x}_T=\bar{x}_T((s_1^{\varepsilon},s_2^{\varepsilon},s_3^{\varepsilon}),x_1)$ and $x_1$ is an arbitrary element of $S$.
\end{Theorem}

Now, we consider the case when two players play
good strategies  and the third one deviates and
chooses an arbitrary strategy. The deviating player does not improve their payoff more then $\frac{2}{3}\varepsilon$, where the positive constant $\varepsilon$ can be chosen arbitrarily small by the two non-deviating players.
\begin{Theorem}\label{T4}
Suppose that the first and the second player choose the $\epsilon$-good
strategies $s_1^\varepsilon$, $s_2^\varepsilon$
and the third player plays an arbitrary strategy $s_3:S\rightarrow
\{I,NI\}$. Then
\begin{equation}
\limsup_{T\to \infty }\bar{x}^3_{T}\leq
p_3+{\epsilon}\frac{2}{3},\label{ne}
\end{equation}
where $\bar{x}_T=\bar{x}_T((s_1^{\varepsilon},s_2^{\varepsilon},s_3),x_1)$ and $x_1$ is an arbitrary element of $S$.
\end{Theorem}
 At the end of the section we
show an example of the third player strategy, such that the
upper limit of his average payoffs is strictly
grater than $p_3$.

Now, we consider the case when two players deviate.
\begin{Theorem}\label{T2}
Suppose  $s_1^\varepsilon$ is the $\epsilon$-good strategy
for the first player and $s_2$, $s_3$ are arbitrary strategies. Then
\begin{equation}
\liminf_{T\rightarrow \infty} \bar{x}^1_T\geq r_0,\label{lEN}
\end{equation}
\begin{equation}
\limsup_{T\rightarrow \infty} (\bar{x}^2_T+\bar{x}^3_T )\leq
2p_3, \label{ngOP}
\end{equation}
\begin{equation}
\lim_{T\to\infty} \mbox{dist}(\bar{x}_T,\,V_1)=0, \label{attr}
\end{equation}
where $\bar{x}_T=\bar{x}_T((s_1^{\varepsilon},s_2,s_3),x_1)$ and $x_1$ is an arbitrary element of $S$.
\end{Theorem}

Suppose that the payoff in the repeated game is defined as the Banach limit of average payoffs.
The inequality (\ref{ngOP}) provides that if two players deviate then at least one of them will not improve his payoff.
 Conclusions (\ref{lEN}) and (\ref{attr}) mean that the good strategy is safe, i.e. the non-deviating player's payoff is not smaller than the  Nash equilibrium payoff in the stage game and, moreover, the deviating player's payoff is not greater than the non-deviating player's payoff plus $\varepsilon$ (comp. Proposition \ref{banachlimit}).

By Theorems \ref{T3} -- \ref{T2}, we obtain
\begin{Corollary}
 The strategy profile $s^\varepsilon=(s_1^{\varepsilon},s_2^{\varepsilon},s_3^{\varepsilon})$ satisfies (\ref{war1}-\ref{war5}). If we define the payoff in the repeated game as a Banach limit of average payoffs, i.e $Lim\bar{x}_T$ then the strategy profile $s^\varepsilon$ is a safe and strong $\varepsilon$ Nash equilibrium.

 \end{Corollary}

Below we provide some elementary properties of sets $V_i$ that are used in the definition of good strategies.
We assume that $i$, $j$, $k$ are pairwise different elements of the set of players $\{1,\,2,\,3\}$. We shall use the following notations
\[
\begin{array}{l}
V^3=\bigcap_{i=1}^3 V_i\\
 V^2_i=(S\setminus V_i)\cap V_j \cap V_k\\
 V^1_i= V_i\cap (S\setminus V_j)\cap (S\setminus V_k)
 \end{array}
 \]
 If each player plays good strategy then
 \begin{equation}\label{fi}
 \varphi(x)= \left\{
\begin{array}{ll}
B &\mbox{ if } x\in V^3 \\
C^1_i & \mbox{ if } x\in V^1_i \quad \mbox{ for } i\in \{1,2,3\}\\
C^2_i & \mbox{ if } x\in V^2_i \quad \mbox{ for } i\in \{1,2,3\}
\end{array}
\right.
 \end{equation}
 \begin{Proposition}\label{Properties1}
 Suppose that player $i$ plays $\varepsilon_i$-good strategy for $i=1,\,2,\,3$. Then
 \begin{equation}\label{P1}
 \Omega_i\cap W_i=\emptyset
 \end{equation}
 \begin{equation}\label{P2}
 V^1_i\subset \Omega_i
 \end{equation}
where
 \[
 \Omega_i=\{x\in S:x_i=\max\{x_1,x_2,x_3\}\}
   \]
  If we  assume that $\varepsilon_i=\varepsilon_j\,(=:\varepsilon)$ then
  \begin{equation}\label{P4}
V_i\cap(S\setminus V_j)\subset \Omega_i\cup\Phi_j
 \end{equation}
If we  assume that $\varepsilon_i=\varepsilon_j=\varepsilon_k\,(=:\varepsilon)$ then for every $i\in\{1,\,2,\,3\}$ we have
 \begin{equation}\label{P3}
V^2_i\subset \Phi_i,
 \end{equation}
 where
\[
  \Phi_i=\{x\in S:x_i=\min\{x_1,x_2,x_3\}\}
  \]
 \end{Proposition}
{\bf Proof : }
If $x_i<r_0$ and $x_i=\max\{x_1,x_2,x_3\}$ then $x_1+x_2+x_3<3r_0$. If $x_j+x_k>2p_3$ and $x_i=\max\{x_1,x_2,x_3\}$ then
$x_1+x_2+x_3>3p_3$. Since $3r_0\leq x_1+x_2+x_3\leq3p_3$ for $x\in S$, we obtain (\ref{P1}).\\
As $(S\setminus V_j)\cap\Omega_j=\emptyset$ and $(S\setminus V_k)\cap\Omega_k=\emptyset$ we have $(S\setminus V_j)\cap (S\setminus V_k)\subset S\setminus (\Omega_j\cup\Omega_k)\subset \Omega_i$, and consequently we obtain (\ref{P2}).\\
To prove  conclusion (\ref{P4}) we take $i=1$, $j=2$. As $(W_2\setminus W_1)\cap\Phi_1=\emptyset$ and $(\Omega_3\setminus \Phi_1)\subset\Phi_2$ we obtain $W_2\setminus W_1\subset (\Omega_1\cup\Omega_3)\setminus\Phi_1\subset\Omega_1\cup\Phi_2$. If $x\in\Omega^\varepsilon_1\setminus\Omega^\varepsilon_2$ then either
\[
x_2\leq x_1-\varepsilon
\]
or
\[
x_2\leq x_3-\varepsilon \mbox{ and } x_1>x_3-\varepsilon\;(x\in\Omega^\varepsilon_1)
\]
In both cases we obtain $x_2<x_1$ and thus $x\in\Omega_1\cup\Phi_2$. Since $V_1\cap(S\setminus V_2)\subset(\Omega^\varepsilon_1\setminus\Omega^\varepsilon_2)\cup(W_2\setminus W_1)$, we conclude that
\[
V_1\cap(S\setminus V_2)\subset \Omega_1\cup\Phi_2.
\]
If $x_i\leq x_j-\varepsilon$ ($x\notin\Omega^{\varepsilon}_i$) and $x_k>x_j-\varepsilon$ ($x\in\Omega^{\varepsilon}_k$) then $x\in \Phi_i$.
If $x_i\leq x_k-\varepsilon$ ($x\notin\Omega^{\varepsilon}_i$) and $x_j>x_k-\varepsilon$ ($x\in\Omega^{\varepsilon}_j$) then $x\in \Phi_i$.
Thus $(S\setminus\Omega^{\varepsilon}_i)\cap\Omega^{\varepsilon}_j\cap\Omega^{\varepsilon}_k\subset\Phi_i$\\
If $x\in V_j$ then $x\notin  W_j$ and hence $x_j\geq r_0$ and $x_i+x_k\leq 2p_3$. If $x\in W_i\cap V_j\cap V_k$ then either $x_i<r_0,\; x_j\geq r_0,\;x_k\geq r_0$ or $ x_j+x_k>2p_3,\;x_i+x_k\leq 2p_3,\; x_i+x_j\leq 2p_3$.
In both cases we deduce that $x\in\Phi_i$. So $V^2_i\subset \Phi_i$.

\begin{flushright} QED \end{flushright}

First we prove Theorem $\ref{T2}$.\\
{\bf Proof : }
The strategy $s^*_1$ is the $\epsilon$-good strategy, so if $\bar{x}^1_T<r_0$ then
$\bar{x}_T=(\bar{x}^1_T,\bar{x}^2_T,\bar{x}^3_T)\in W_1$ and $s^*_1(\bar{x}_T)=NI$.
It means that the next vector payoff $x_{T+1}$ belongs to the set
$\{A,C^1_2,C^1_3,C^2_1\}$, so $x^1_{T+1}\in \{r_0,r_1,r_2\}$, i.e.
$x^1_{T+1}\geq r_0$ (see $(\ref{nierpodst})$). By Corollary \ref{lBC} we obtain that
$\limsup_{T\rightarrow \infty} -\bar{x}^1_T\leq  -r_0$, so
$\liminf_{T\rightarrow \infty} \bar{x}^1_T\geq r_0$.

Similarly, if $\bar{x}^2_T+\bar{x}^3_T>2p_3$ then $s^*_1(\bar{x}_T)=NI$. Thus the sum $x^2_{T+1}+x^3_{T+1}$ is one of the numbers: $2r_0,p_1+r_1,2p_2$. From the assumptions $(\ref{nierpodst})$, $(\ref{nierinvest})$ and  $(\ref{nierstrong})$, it follows that $x^2_{T+1}+x^3_{T+1}\leq 2p_3$.  By Corollary
\ref{lBC}, we get
\begin{displaymath}
\limsup_{T\rightarrow \infty} \bar{x}^2_T+\bar{x}^3_T \leq
2p_3.
\end{displaymath}
If $x\in S\setminus V_1$ then $s^*_1(\bar{x}_T)=NI$. So, $\varphi(x)=G((s^*_1,s_2,s_3)(x))\in \{A,C^1_2,C^1_3,C^2_1\}\subset V_1$. By Corollary \ref{2bc}, the set $V_1$ is a weak attractor for $\varphi$.

\begin{flushright} QED \end{flushright}

Let $\pi_u:\mathbb{R}^3\rightarrow u$ be the orthogonal  projection  onto the line $u=\{x\in
\mathbb{R}^3: x_1=x_2=x_3\}$ and $\pi_P:\mathbb{R}^3\rightarrow P$ be the orthogonal
projection  onto the plane $P=\{x\in \mathbb{R}^3:x_1+x_2+x_3=0\}$. Obviously
$\pi_u(x)=\left(\frac{x_1+x_2+x_3}{3},
\frac{x_1+x_2+x_3}{3},\frac{x_1+x_2+x_3}{3}\right)$ and
$
\pi_P(x)=x-\pi_u(x)
$.
In the remainder of the section we denote the projection of a point (a set) $A$ onto the plane $P$ by $\tilde{A}$, i.e. $\tilde{A}=\pi_P(A)$.
The projection of the set $S$ onto the plane  $P$:
\begin{displaymath}
\tilde{S}:=\pi_P(S)
\end{displaymath}
is the convex hull of the hexagon with successive vertexes $\tilde{C^1_1}$, $\tilde{C^2_2}$, $\tilde{C^1_3}$, $\tilde{C^2_1}$, $\tilde{C^1_2}$, $\tilde{C^2_3}$.

 Set
\begin{displaymath}
\begin{array}{lll}
v_1=\frac{1}{\sqrt{2}}(0,-1,1) & \qquad & v_4=-v_1 \\
v_2=\frac{1}{\sqrt{2}}(-1,0,1) & \qquad & v_5=-v_2 \\
v_3=\frac{1}{\sqrt{2}}(-1,1,0) & \qquad & v_6=-v_3
\end{array}
\end{displaymath}
and
\begin{displaymath}
\Delta_c(K)=\bigcap_{i\in K} \{y\in \tilde{S}:\langle v_i,y\rangle <c\}
\end{displaymath}
where $K\subset\{1,\ldots,6\}$ and $c>0$. One can easy check that
\[
\begin{array}{l}
x\in\Omega^\varepsilon_1\Leftrightarrow\pi_P(x)\in\Delta_c(\{2,3\})\\
x\in\Omega^\varepsilon_2\Leftrightarrow\pi_P(x)\in\Delta_c(\{1,6\})\\
x\in\Omega^\varepsilon_3\Leftrightarrow\pi_P(x)\in\Delta_c(\{4,5\})
\end{array}
\]
where $c=\frac{\varepsilon}{\sqrt{2}}$. Setting
$
\Omega^{\epsilon}=\bigcap_{i=1}^3\Omega_i^{\epsilon}
$
and $\Delta_c=\Delta_c(\{1,\ldots,6\})$ we obtain

\begin{equation}
x\in
\Omega^{\epsilon} \Leftrightarrow \pi_P(x)\in \Delta_c.\label{rownowa}
\end{equation}

Now, we are able to prove Theorem $\ref{T3}$.\\
{\bf Proof : }
Fix $x_1\in S$. It is  sufficient to show that in the sequence $\bar{x}_T=\bar{x}_T(s^*,x_1)$ there exists an element $\bar{x}_N$ belonging to $V^3$, where $s^*=(s_1^{\varepsilon},s_2^{\varepsilon},s_3^{\varepsilon})$.  Indeed, if $\bar{x}_N\in V^3$ then $\bar{x}_{N+k}=\frac{N}{N+k}\bar{x}_N+\frac{k}{N+k}B$, so $\lim_{T\to\infty}\bar{x}_T=B$.\\
First we show that almost all elements of the sequence $\bar{x}_T$ belong to $\Omega^\eta=\bigcap_{i=1}^3\Omega_i^{\eta}$, for every $\eta>0$.\\
The map $\varphi$ given by (\ref{fi}) is determined by the strategy profile $s^*$, i.e. $\varphi=G\circ s^*$.
Consider $\tilde{\varphi}:\tilde{S}\multimap \tilde{S}$ and $V:P\to\mathbb{R}$ given by
\[
\begin{array}{l}
\tilde{\varphi}(x)=\{\pi_P(\varphi(y)):\;\pi_P(y)=x\}\\
V(x)=\max\left\{\langle v_i,x\rangle : i=1,\dots,6\right\}
\end{array}
\]
We verify that $V$ is a Lapunov type function for $\tilde{\varphi}$ with the constant $c$, for an arbitrary $c>0$. Let us fix $x\in\tilde{S}$ such that $\langle v_1,x\rangle >0$. If $y\in S$ and $\pi_P(y)=x$ then $\langle v_1,y\rangle=\langle v_1,x\rangle$.Thus $y_3-y_2>0$ and therefore $y\notin\Omega_2\cup\Phi_3$. By (\ref{P2}), (\ref{P3}), we have $y\notin V^1_2\cup V^2_3$. Since $\varphi(y)\in\{C^1_1,\,C^1_3,\,C^2_1,\,C^2_2,\,B\}$, we obtain $\langle v_1,\omega\rangle\leq 0$ for $\omega\in\tilde{\varphi}(x))$. We use similar arguments to show that  if $\langle v_i,x\rangle> 0$ and $\omega\in\tilde{\varphi}(x)$ then $\langle v_i,\omega\rangle\leq 0$, for $i=2,\ldots,6$.\\
Fix $\eta<\min\{\epsilon,\,p_3-\frac{2p_2+r_2}{3},\,\frac{p_1+2r_1}{3}-r_0\}$ .
By Proposition \ref{lfm} and (\ref{rownowa}), there exists $N$ such that $\bar{x}_n\in\Omega^\eta$ for $n>N$. We claim that there exists $M>N$ such that $\bar{x}_M\in V^3$. Suppose to the contrary that $\bar{x}_M\notin V^3$ for every $M>N$. Then $\varphi(\bar{x}_M)\in\{C^1_1,\,C^1_2,\,C^1_3,\,C^2_1,\,C^2_2,\,C^2_3\}$ for $M>N$. By Proposition \ref{lA}, we obtain that $z_M\in (\frac{p_1+2r_1-\eta}{3},\,\frac{2p_2+r_2+\eta}{3})$ for $M$ sufficiently large, where the point $(z_M,z_M,z_M)$ is the projection of $\bar{x}_M$ onto $u$. \\
But, if $x\in\Omega^\eta\setminus V^3$ then $\frac{x_1+x_2+x_3}{3}\notin(\frac{p_1+2r_1-\eta}{3},\,\frac{2p_2+r_2+\eta}{3})$.
Indeed, if $x_i+x_j>2p_3$ and $x\in\Omega^\eta$ then $\frac{x_1+x_2+x_3}{3}>p_3-\frac{\eta}{3}$. If $x_i<r_0$ and $x\in\Omega^\eta$ then $\frac{x_1+x_2+x_3}{3}<r_0+\frac{2}{3}\eta$.
\begin{flushright} QED \end{flushright}

Now, we are in a position to prove Theorem $\ref{T4}$.\\
{\bf Proof : }
Let $x_1\in S$ and $\eta>0$. Our aim is to prove that almost all
elements of the sequence
$\bar{x}_T=\bar{x}_T((s^\varepsilon_1,\,s^\varepsilon_2,s_3),\;x_1)$
belongs to
$\Omega^{\varepsilon+\eta}:=\Omega^{\varepsilon+\eta}_1\cap\Omega^{\varepsilon+\eta}_2$.
We have
\begin{equation}\label{od}
x\in\Omega^{\varepsilon+\eta}\Leftrightarrow
\pi_P(x)\in\Delta_c(\{1,\,2,\,3,\,6\}),
\end{equation}
where $c:=\frac{\varepsilon+\eta}{\sqrt{2}}$. We show that the
function $V^*:P\rightarrow \mathbb{R}$ given by
\[
V^*(x)=\max\{\langle v_i,x\rangle : i=1,2,3,6\}
\]
is the Lapunov type function for
$\tilde{\varphi^*}:\tilde{S}\multimap\tilde{S}$ with the constant
$c$, where
\[
\tilde{\varphi^*}=\{\pi_P(z);\;z\in \varphi^*(y),\;\pi_P(y)=x\}
\]
and
\begin{displaymath}
\varphi^*(y)=\left\{
\begin{array}{ll}
\{C^1_1,C^2_2\} &\mbox{ if } y\in V_1\cap (S\setminus V_2) \\
\quad & \quad \\
\{C^1_2,C^2_1\} &\mbox{ if } y\in (S\setminus V_1)\cap V_2 \\
\quad & \quad \\
\{B,C^2_3\} &\mbox{ if } y\in V_1\cap V_2 \\
\quad & \quad \\
\{A,C^1_3\} &\mbox{ if } y\in (S\setminus V_1) \cap (S\setminus V_2) \\
\end{array}\right..
\end{displaymath}
The map $\varphi:S\rightarrow S$ induced by the profile
$(s^\varepsilon_1,\,s^\varepsilon_2,s_3)$ is a selection of
$\varphi^*$.

If $\langle v_6,x\rangle>0$ ($x\in\tilde{S}$) and $\pi_P(y)=x$
($y\in S$) then $y_1>y_2$ and thus $y\notin \Omega_2\cup\Phi_1$.
By (\ref{P4}), we have $V_2\cap(S\setminus V_1)\subset
\Omega_2\cup\Phi_1$. Thus
$\varphi^*(y)\cap\{C^1_2,\,C^2_1\}=\emptyset$. So, we have
$\langle v_6, \omega\rangle<0$ for $\omega\in
\tilde{\varphi^*}(x)$.\\
Using similar arguments we show that if $\langle v_3,x\rangle>0$
then $\langle v_3,\omega\rangle\leq 0$ for
$\omega\in\tilde{\varphi^*}(x)$.\\
Suppose that $\langle v_1,\,x\rangle\geq V^*(x)-\delta$
($x\in\tilde{S}$) and $\pi(y)=x$ ($y\in S$), where
$\delta<\frac{\eta}{\sqrt{2}}$. Then $\langle v_1 ,x
\rangle=\langle v_1,y \rangle>\frac{\varepsilon}{\sqrt{2}}$. If
$z\in\Omega^\varepsilon_2$ then $\langle v_1,z
\rangle\geq\frac{1}{\sqrt{2}}(z_3-z_2)>\frac{\varepsilon}{\sqrt{2}}$.
Thus, we have $y\notin \Omega^\varepsilon_2\supset V_2$ and
therefore $\varphi^*(y)\subset \{C^1_1,\,C^2_2,\,C^1_3,A\}$. So,
$\langle v_1,\omega \rangle\leq 0$ for $\omega\in
\tilde{\varphi^*|}(x)$.\\
In the similar way we prove that if $\langle v_2,x \rangle\geq
V^*(x)-\delta$ and $\omega\in\tilde{\varphi^*}(x)$ then $\langle
v_2,\omega \rangle\leq 0$.

By Proposition \ref{lfm}, we obtain that almost all elements of the
sequence ($\pi_P(\bar{x}_T)$) belongs to
$\Delta_c(\{1,\,2,\,3,\,6\})$. By (\ref{od}), we have that almost
all elements of the sequence ($\bar{x}_T$) belongs to
$\Omega^{\varepsilon+\eta}$. If $x\in \Omega^{\varepsilon+\eta}$
then $x_1>x_3-(\varepsilon+\eta)$ and $x_2>x_3-(\varepsilon+\eta)$
and so $x_3<p_3+\frac{2}{3}(\varepsilon+\eta)\;$ ($x_1+x_2+x_3\leq 3
p_3$ for $x\in S$).
\begin{flushright} QED \end{flushright}

 {\bf Remark.} Reasoning as in the proofs of Theorem \ref{T4} and Theorem \ref{T2}, we can conclude that good strategies are safe and strong Nash equilibria not only in the class of Smale's strategies, but also if "loyal" players adopt good strategies, then "disloyal" players can  even play the random choice in each repetition. It does not change the properties (\ref{ne}), (\ref{lEN}), (\ref{ngOP}) and(\ref{attr}).

{\bf Example 2}
Let the stage game $G$ be given by:
\begin{displaymath}
\begin{array}{c|c|c}
n & P_I(n) & P_{NI}(n) \\ \hline
0 & - & 20 \\
1 & 10 & 28 \\
2 & 18 & 36 \\
3 & 26 & -
\end{array}.
\end{displaymath}
This game satisfies  conditions $(\ref{nierpodst})$ -- $(\ref{nierantypod})$.

Let $s_i^*:S\rightarrow \{I,NI\}$ be the
$\epsilon$-good strategy for the $i$-th player, $i=1,2$, and
$0<\epsilon<\frac{1}{2}$. Let $Z= conv\{A,B,C^1_3,C^2_3\}=\{x\in
S:x_1=x_2\}$ and $D=(26-\frac{\epsilon}{2},
26-\frac{\epsilon}{2}, 26+\frac{\epsilon}{2})$. We present
the construction of the  third player strategy $s_3^*:S\to\{I,NI\}$   such that
\begin{displaymath}
\lim_{T\to\infty} x_T(s^*,\,x_1)=D\mbox{ for every } x_1\in Z
\end{displaymath}
We have $V_1\cap Z=V_2\cap Z$. We set
\begin{displaymath}
s_3^*(x)=\left\{
\begin{array}{ll}
NI &\mbox{ if } x\in V_1\cap Z\cap co\{B,\,D,\,C^1_3\} \\
I & \mbox{ elsewhere. }
\end{array}
\right.
\end{displaymath}
The map $\varphi$ induced by the strategy profile $s^*=(s_1^*,s_2^*,s_3^*)$ is given by
\begin{displaymath}
\varphi(x)=\left\{
\begin{array}{ll}
B & \mbox{ if } x\in V_1\cap Z\setminus co\{B,\,D,\,C^1_3\} \\
C^2_3 &\mbox{ if } x\in V_1\cap Z\cap co\{B,\,D,\,C^1_3\}\\
C^1_3 & \mbox{ if } x\in Z\setminus V_1
\end{array}
\right..
\end{displaymath}
The values of the map $\varphi$ outside the set $Z$  have no influence onto the trajectory $\bar{x}_T(s^*,\,x_1)$ if $x_1\in Z$.
The map $\varphi:Z\to Z$ satisfies the Blackwell condition for the triangle $co\{C^1_3,\, C^2_3,\,D\}$ in the domain $Z$. The map $\varphi:Z\to Z$ satisfies the Blackwell condition for the sum of intervals $\overline{BD}\cup \overline{D\,C^1_3}$ in the domain $Z$. By Proposition \ref{bt}, the sets
$co\{C^1_3,\, C^2_3,\,D\}$ and $\overline{BD}\cup \overline{D\,C^1_3}$ are week attractors. To conclude that the interval $\overline{BD}$ is a weak attractor we apply Proposition \ref{atractor2} taking $A=\overline{BD}\cup \overline{D\,C^1_3}$ and $B=\overline{BD}$. By Proposition \ref{atractor1}, the intersection of weak attractors $co\{C^1_3,\, C^2_3,\,D\}$ and $A=\overline{BD}$ is a weak attractor. The intersection equals to the  set $\{D\}$.

\section{Concluding remarks}\label{sec:4}

This paper is concerned with the specific model of social dilemmas. Such models have a very special place in game theory as they describe real social problems of modern world: resources depletion, pollution and overpopulation. The main characteristic of such models is that each player gain more by not cooperating when opponents fix their choices and all individuals are better off if all cooperate. The lack of optimality of Nash equilibrium is the most interesting problem, because as we can observe in the real world, people are keen to cooperate with each other on the certain conditions. Axelrode showed in \cite{ax} that strategies that effectively encourage people to cooperate are: nice, forgiving, retaliatory and are found on simple rules.

The key idea in our approach is to apply Smale's idea for 3-payer extension of Presoner's Dilemma. Our strategies are deterministic and satisfy conditions that are postulated by Axelrode in \cite{ax}. What is more, $\epsilon$-good strategies satisfy condition (\ref{war5}) which guarantee that using this strategy our payoff shall not be different than our opponents payoffs for more than $\epsilon$. This constant  $\epsilon$ is totally controlled by the player who choose it. This property is not received by any other author.

Our future aim is to extend the idea presented in \cite{pz} onto  the type of games considered in the paper - three players repeated social dilemmas. The idea is as follows. We would like to analyze  the repeated three players game by evolutionary games methods. To achieve this goal,   we threat the repeated game as a new game in which a player action is a  point in the $\beta$-core of the original game.  Using methods presented in the paper each point from the $\beta$-core should determine $\varepsilon$-good strategy. The main difficulty is to obtain the payoff in the case when players choose different points in the $\beta$-core. The payoffs in the new game are determined by the payoff in the repeated game.

\appendix
\section{Appendix}\label{sec:A}

In this Appendix we shall present the proof of Proposition \ref{lfm}. We start with the necessary theorem.

\begin{Theorem}\label{T1}
If $S$ is a bounded convex subset of $H$ and
$V:H\to\mathbb{R}$ given by (\ref{function V}) is the
Lapunov type function for $\varphi:S\multimap S$ with the constant
$c>0$, then
\[
\exists \gamma>0,\quad \exists \alpha_0 >0, \quad \forall \alpha \in
[0,\alpha_0] ,\quad \forall x\in S\setminus \Delta_c ,\quad \forall
\omega \in \varphi (x)
\]
\begin{equation}
V(\alpha \omega + (1-\alpha)x) \leq V(x) -\alpha
\gamma.\nonumber
\end{equation}
\end{Theorem}

{\bf Proof : }
 By $(\ref{Lapunov})$ we choose $\delta\in(0,\,c)$. Let $M=\sup\{|x|:\,x\in S\}$. For $x\in
S\setminus \Delta_c$ we define a set of indexes $I(x)$ by
\[
I(x)=\{j\in \{1,...,k\}: V_j(x)\geq V(x) -\delta\},
\]
and a subset $O_i$ of S, related to the fixed index $i$:
\begin{equation}
O_i:=\{x\in S\setminus \Delta_c: V_i(x) \geq V(x) -
\delta\}.\label{zbioryO}
\end{equation}
If  $x\in O_i$ then $V_i(x)>0$ and  $\langle p_i, \omega \rangle
\leq 0$ for all $\omega \in \varphi(x)$. Obviously, $i\in
I(x)$ is equivalent to $ x\in O_i $ for $x\in S\setminus\Delta_c$.

We fix positive constants: $r$, $\gamma$ and $\alpha_0$ such that
\[
r<\frac{\delta}{2}, \qquad \gamma<c-\delta, \qquad \alpha_0<
\min\left\{\frac{c-\delta - \gamma}{c-\delta + M}, \frac{r}{2M},1\right\}
\]
and take an arbitrary $x\in S\setminus \Delta_c$. The following
condition holds true
\begin{equation}
\forall y\in B(x,r)=\{y\in S\setminus \Delta_c:||x-y||<r\} \quad
\exists i\in I(x) \quad V(y)=V_i(y). \label{row1}
\end{equation}
Indeed, if $j\notin I(x)$, then $V_j(x)<V(x)-\delta$. Since $V_j$
and $V$ are   lipschitz continuous  with the constant L=1, we get
$V_j(y) <V(y)$. Therefore, there exists $i\in I(x)$ such that
$V(y)=V_i(y)$.

If $i\in I(x)$ then $V_i(x)\geq V(x)-\delta\geq
c-\delta$ and $\langle p_i,\omega\rangle\leq 0$ for
$\omega\in\varphi(x)$. Let $\alpha\in[0,\alpha_0]$ then
$x_\alpha:=\alpha\omega+(1-\alpha)x\in B(x,r)$ and
$V_i(x_{\alpha_0})\leq V_i(x_\alpha)$. Moreover,
\begin{displaymath}
V_i(x_{\alpha})  \geq  V_i(x_{\alpha_0})
 \geq  -\alpha_0||\omega|| +( 1-\alpha_0)(c-\delta)
\geq c-\delta -\alpha_0(c-\delta +M)
 \geq \gamma\nonumber
\end{displaymath}
and
\begin{displaymath}
V_i(x_{\alpha})\leq (1-\alpha)V_i(x)\leq V_i(x)-\alpha
\gamma.
\end{displaymath}
Thus we have obtained that
\begin{equation}
\forall i\in I(x),\quad \forall \alpha \in [0,\alpha_0], \quad
\forall \omega \in \varphi(x), \quad V_i(x_{\alpha})\leq
V_i(x) - \alpha \gamma.\label{row4}
\end{equation}
The function $V$ has the following property: if $V(a)=V_i(a)$ and
$V(b)=V_i(b)$ then $V(\lambda a+(1-\lambda)b)=V_i(\lambda
a+(1-\lambda)b)$ for $\lambda\in[0,1]$ so the set
\[
\{\alpha \in [0,\alpha_0] : V_i(\alpha \omega +(1-\alpha)x)=
V(\alpha \omega +(1-\alpha)x)\}.
\]
is a closed segment. By (\ref{row1}) there exists $s\leq k$
and a partition $0=\beta_0<\beta_1<...<\beta_s=\alpha_0$ such that
\begin{equation}
\forall j\in \{0,...,s-1\}, \, \exists i=i(j)\in I(x), \,
\forall\alpha \in [\beta_j,\beta_{j+1}], \quad
V(x_{\alpha})=V_i(x_{\alpha}).\label{ll4}
\end{equation}
Let $\alpha \in [\beta_0,\beta_1]$. In view of $(\ref{ll4})$ there
exists $i=i(0)\in I(x)$ such that $V(x_{\alpha})=V_i(x_{\alpha})$
and by $(\ref{row4})$:
\[
V(x_\alpha)=V_i(x_{\alpha})\leq V_i(x)-\alpha \gamma = V(x)
-\alpha \gamma.
\]
Suppose that
\[
V(x_{\alpha})\leq V(x) -\alpha \gamma, \quad
\forall k=1,...,j-1 \quad \forall \alpha \in
[\beta_k,\beta_{k+1}]
\]
and take $\alpha \in [\beta_j,\beta_{j+1}]$. By $(\ref{ll4})$
there exists $i=i(j)$ such that $V(x_{\alpha})=V_i(x_{\alpha})$.
Since $\alpha \in [\beta_j,\beta_{j+1}]$, there exists $\xi\in
[0,1]$ such that $x_{\alpha}=\xi \omega +(1-\xi)x_{\beta_j}$.
Therefore,
\[
x_{\alpha}=\xi \omega + (1-\xi)(\beta_j \omega +(1-\beta_j)x)=
(\xi+(1-\xi)\beta_j)\omega + (1-\xi)(1-\beta_j)x,
\]
so
\[
\alpha= \xi +(1-\xi)\beta_j \leq \xi + \beta_j.
\]
It is obvious that
\[
V(x_{\alpha})=V_i(x_{\alpha})=V_i(\xi \omega +(1-\xi)x_{\beta_j})
= \xi \left<\omega, p_i\right>+(1-\xi)\left<x_{\beta_j},p_i\right>
\]
\[
\leq (1-\xi)V_i(x_{\beta_j})= V_i(x_{\beta_j})-\xi
V_i(x_{\beta_j}) \leq V_i(x_{\beta_j}) -\xi \gamma.
\]
Then we get
\[
V_i(x_{\beta_j}) -\xi \gamma \leq V(x)-\beta_j\gamma  -\xi
\gamma = V(x) - \gamma (\beta_j+\xi) \leq V(x) -\alpha
\gamma,
\]
hence,
\[
V(\alpha \omega +(1-\alpha)x) \leq V(x) -\alpha \gamma.
\]
\begin{flushright} QED \end{flushright}

The proof of Proposition \ref{lfm}.

{\bf Proof : }
Fix $(\bar{x}_n)_{n=1}^\infty$ satisfying (\ref{zm}).  First, we
prove that
\begin{equation}\label{w1}
\forall M,\,\exists N\geq M,\;\;\bar{x}_N\in\Delta_c
\end{equation}
Suppose, contrary to our claim, that $\bar{x}_n\notin\Delta_c$ for
$n\geq m$. We choose $k\geq m$ such that $\frac{1}{k}<\alpha_0$,
where $\alpha_0$ and $\gamma$ are given by Theorem \ref{T1}. Thus
\[
V(\bar{x}_{k+l+1})=V\left(\frac{1}{k+l+1}x_{k+l+1}+\frac{k+l}{k+l+1}\bar{x}_{k+l}\right)\leq
V(\bar{x}_{k+l})-\gamma\frac{1}{k+l+1}\leq
\]
\[ \ldots\leq V(\bar{x}_k)-\gamma\left(\frac{1}{k+l+1}+\ldots+\frac{1}{k+1}\right) \xrightarrow[\l\to\infty]{}-\infty
\]
which contradicts to the assumption that $V(\bar{x}_n)\geq c$ for
$n\geq m$.

Fix $c_1>c$. By Proposition \ref{l1}, we choose $M$ such that
$|\bar{x}_{l+1}-\bar{x}_l|<c_1-c$ for $l\geq M$. By (\ref{w1}),
there exists $N\geq M$ such that $\bar{x}_N\in\Delta_c$. If
$\bar{x}_{N+l}\in\Delta_c$ then $V(\bar{x}_{N+l+1})\leq
V(\bar{x}_{N+l})+|\bar{x}_{N+l+1}-\bar{x}_{N+l}|<c_1$. If
$\bar{x}_{N+l}\in\Delta_{c_1}\setminus\Delta_c$ then
$V(\bar{x}_{N+l+1})\leq V(\bar{x}_{N+l})<c_1$.

\begin{flushright} QED \end{flushright}

\end{document}